\documentstyle[12pt]{amsart}

\newtheorem{conjecture}{Conjecture}[section]
\newtheorem{theorem}{Theorem}[section]
\newtheorem{example}{Example}[section]
\newtheorem{definition}{Definition}[section]

\newtheorem{lemma}{Lemma}[section]
\newtheorem{corollary}{Corollary}[section]
\newenvironment{proof}{\noindent {\bf {Proof:}}}{$\Box$}
\newenvironment{proof-5-2}{\noindent {\bf {Proof of Theorem
5.2:}}}{$\Box$}
\newenvironment{proof-5-1}{\noindent {\bf {Proof of Theorem
5.1:}}}{$\Box$}
\newenvironment{proof-7-2}{\noindent {\bf {Proof of Theorem
7.2:}}}{$\Box$}
\newenvironment{proof-7-1}{\noindent {\bf {Proof of Theorem
7.1:}}}{$\Box$}
\newcommand{\T}{(t^{\frac12} + t^{-\frac12})}
\newcommand{\Q}{(q-1)}
\newcommand{\hr}{(\frac{r-1}{2})}
\renewcommand{\r}{{\text{\em Res}}}
\numberwithin{equation}{section}

\begin{document}

\title[On Ohtsuki's Invariants]
{On Ohtsuki's Invariants of \\ Integral Homology 3-Spheres, I}
\author{Xiao-Song Lin}
\address{Department of Mathematics, University of California, Riverside,
CA 92521}
\email{xl@@math.ucr.edu}
\author {Zhenghan Wang}
\address {Department of Mathematics, University of Michigan, Ann Arbor,
MI 48109}
\email{Zhenghan.Wang@@math.lsa.umich.edu}
\date{September, 1995. Revised: April, 1996}
\thanks{The first author is supported in part by NSF and
the second author is supported
by a Rackham Faculty Summer Fellowship at University of Michigan.}
\begin{abstract} An attempt is made to conceptualize the derivation as well
as to facilitate the
computation of Ohtsuki's
rational invariants $\lambda_n$ of integral homology 3-spheres extracted
from Reshetikhin-Turaev $SU(2)$ quantum invariants.
Several interesting consequences will follow from our computation
of $\lambda_2$. One of them says that $\lambda_2$ is always an integer
divisible by 3. It seems interesting to compare this result
with the fact shown by
Murakami that $\lambda_{1}$ is 6 times the Casson invariant. Other
consequences include some general criteria for
distinguishing homology 3-spheres
obtained from surgery on knots by using the Jones polynomial.
\end{abstract}
\maketitle

\section{Introduction}

In [O1], Ohtsuki extracted a series of rational
topological invariants of oriented homology 3-spheres
from the $SU(2)$ quantum invariants of Reshetikhin and Turaev [RT].
Physically, they correspond to the coefficients of the
asymptotic expansion of Witten's Chern-Simons path integral at
the trivial connection as shown by Rozansky [R1,2].
The original derivation of these invariants, though,
seems too complicated to make any practical computation possible.
We will present here an attempt to conceptualize as well
as to simplify the derivation of Ohtsuki's invariants.
The simplified derivation will lead to formulae for Ohtsuki's invariants,
which are very practical and straightforward when it comes to
actual computation. Various interesting consequences will
follow from our formulae for Ohtsuki's invariants.

To be more precise, Ohtsuki [O1] defined a formal power series
$$\tau(M)=1+\sum_{n=1}^\infty\lambda_n(t-1)^n$$
with rational coefficients for every oriented integral homology 3-sphere $M$.
It is, in our terminology, the {\it Fermat limit} of Reshetikhin-Turaev
$SU(2)$
quantum invariants $\{\tau_r(M)\}$ of $M$ (See the definition in Section 3).
The notion of Fermat limits has appeared implicitly in many places
in the literature on quantum 3-manifold invariants.
But it seems that Murakami was the first to use it
effectively in identifying $\lambda_1$ with
6 times the Casson invariant [Mu].

Ohtsuki's invariants $\{\lambda_n\}$ are closely related
with the developing theory of finite type 3-manifold invariants
[Ha, G, GLe1, GLe2, GLi, GO, O2].
As we will see in the sequel to this paper,
some special properties of the colored Jones polynomial
of algebraically split links, i.e., links with all pairwise
linking numbers zero, play a key role in this relationship.

One of the consequences of our computation of Ohtsuki's invariants
concerns their integrality. In particular, we have

\begin{theorem} $\lambda_2(M)\in 3\Bbb Z$ for all integral
homology 3-spheres $M$.
\end{theorem}

This theorem, together with the fact that $\lambda_1\in6\Bbb Z$, suggests the following conjecture.

\begin{conjecture} $n!\cdot \lambda_n(M) \in 6\Bbb Z$ for all
integral homology 3-spheres $M$.
\end{conjecture}

See Section 4 for a relevant conjecture about the integrality of coefficients
of a certain variation of the Jones polynomial on algebraically split links
(Conjecture 4.1). We will say more about this conjecture in the course of our
discussion.

The Casson invariant for homology 3-spheres comes originally from an
algebraic counting of conjugacy classes of irreducible representations of the
fundamental group into $SU(2)$ [AM]. This leads to a general criterion for
detecting knots with property $P$ using the Alexander or Jones polynomial.  To
explain this in some detail, we denote by $J(K,t)$ the 1-variable Jones
polynomial  of a knot $K$ in $S^3$. We write
$$J(K,e^h)=1+\sum_{n=2}^\infty (-1)^n\frac{v_n(K)}{n!}h^n.$$
We will also denote by $S^3_{K,p/q}$ the 3-manifold obtained from
Dehn surgery on the knot $K$ with the surgery coefficient $p/q$.
Then, using the
Casson invariant, we may conclude that if $v_2(K)\neq0$, $S^3_{K,p/q}$
will never be a homotopy 3-sphere unless $p/q=1/0$. We have
several similar criteria for distinguishing homology 3-spheres
obtained from surgery on knots by using other coefficients of the Jones
polynomial in relation with Ohtsuki's invariant
$\lambda_2$. For example, we have the following theorem.

\begin{theorem} Let $K^*$ be the mirror image of $K$. If $v_3(K)\neq0$,
then $S^3_{K,1/n}\neq S^3_{K^*,1/n}$ as unoriented manifolds unless $n=0$.
\end{theorem}

See Corollary 5.3. Needless to say, we don't know what kind of geometrical or
topological obstruction the invariant $\lambda_2$ represents which prevents
these two 3-manifolds $S^3_{K,1/n}$ and $S^3_{K^*,1/n}$ from being
homeomorphic to each other for $n\neq0$.

\medskip

\noindent{\it Acknowledgments:} We would like to express our appreciation
to Shouwu Zhang
for some stimulating discussions about the
notion of Fermat functions, and to Mike Freedman, Yongwu Rong and Ying-Qing Wu
for sharing with us their knowledge about surgery on knots.

\section{Fermat functions and their residues}

Fermat's little theorem tells us that for an integer $a$ and
a prime $r$ such that $a\neq 0\; mod \; r$, $a^{r-1}=1$ mod $r$.
Let us view this theorem from a different angle.

Suppose we have a function $f$ which sends each prime $r$ to $f(r)\in
{\Bbb Z}_r $, where ${\Bbb Z}_{r}$ is the ring of $r$-adic integers.
Suppose $a\in {\Bbb Z}_{r}$ be an r-adic integer, then $\bar{a}$ denotes
the first digit in the r-adic expansion of $a$.
We will call $f(r)$ a {\it Fermat function} if there is a
rational number $\lambda$ independent of $r$ such that
\begin{equation} f(r) \equiv \lambda \qquad\text{mod $r$}
\end{equation}
i.e. $\overline{f(r)}=\overline{\lambda}$
for all sufficiently large primes $r$.
Here we should notice that for the fixed rational
number $\lambda$, we have $\lambda\in {\Bbb Z}_{r}$
if $r$ is sufficiently large.
So (2.1) makes sense for sufficiently large $r$.

If $f$ is a Fermat function, we will call the rational number
$\lambda$ the {\it residue} of $f$ and it is denoted by $\text{Res}(f)$.

\begin{lemma} The residue of a Fermat function is unique.
\end{lemma}

\begin{proof} If there is another rational number $\lambda'$ such that
$f(r)\equiv\lambda'$ mod $r$ for all
 sufficiently large primes $r$, then the numerator of
the rational number $\lambda-\lambda'$ will be divisible
by all sufficiently large primes $r$.
Therefore, $\lambda-\lambda'=0$,  i.e. $\lambda=\lambda'$.
\end{proof}

With these terminologies, Fermat's Little Theorem can be rephrased by saying
that if $a\neq 0\; mod \; r$, then
the function $f(r)=a^{r-1}$ is a Fermat function with $\text{Res}(f)=1$.
There are certainly plenty of other Fermat functions. For example, the function
$f(r)=(r-1)/2$
is a Fermat function whose residue is $-1/2$. Also, the function $f(r)=(r-1)!$ is a
Fermat function because of Wilson's theorem, which says that
$(r-1)!\equiv-1$ mod $r$ [HW]. On the other hand, an example of a
non-Fermat function is given by $f(r)=(\frac{r-1}{2})!$.  See [HW] for a
discussion of the residue of $(\frac{r-1}{2})!$, which turns out to be
dependent on $r$.

The following lemma comes directly from the definition.

\begin{lemma} Suppose that $f$ and $g$ are Fermat functions. Then
\begin{enumerate}
\item for rational numbers $\alpha$ and $\beta$, $\alpha f+\beta g$ is a
Fermat function whose residue is $\alpha \r (f)+\beta \r (g)$;
\item $f\cdot g$ is a Fermat function whose residue is $\r (f)\r (g)$;
\item if $\r (g)\neq 0$, then $f/g$ is a Fermat function whose residue is
$\r (f)/\r (g)$.
\end{enumerate}
\end{lemma}

Thus, in particular, every polynomial function of $r$ with rational
coefficients is a Fermat function whose residue is its constant term.
And every rational function with rational coefficients
is a Fermat function if the constant term of the denominator
is not zero. Its residue will be the value of the function at $r=0$.

The following examples of Fermat functions and the
calculation of their residues
are very important to our discussion of Ohtsuki's invariants.
So we single them out first.
We will use the notation
$${\Bbb Z}_{(m)}=\begin{cases}
{\Bbb Z}[\frac12,\frac13,\dots,\frac1{m-1}]& \text{for $m\geq3$,}\\
{\Bbb Z} & \text{for $m\leq2$}
\end{cases}
$$
with $m$ a positive integer.

\noindent{\bf Example 2.1.} For a fixed integer $k$, the function
\begin{equation}
D_k(r)=\frac{\left(\frac{r-1}2\right)!}{\left(\frac{r-1}2-k\right)!}
\end{equation}
is a Fermat function whose residue is given by
$$\text{Res}(D_k)=\begin{cases}
(-\frac12)(-\frac12-1) \cdots (-\frac12-(k-1))\in{\Bbb Z}_{(3)}& \text{for $k>0$,}\\
1  & \text{for $k=0$,}\\
\frac{1}{(-\frac12+1)\cdots (-\frac12-k)}\in {\Bbb Z}_{(-2k)} &\text{for $k<0$.}
\end{cases}$$

\medskip

\noindent{\bf Example 2.2.} We denote by $\sigma_j(x_1,\dots,x_n)$  ($1\leq
j\leq n$) the $j$-th elementary symmetric polynomial in the indeterminates
$x_1,\dots,x_n$, and $\sigma_0(x_1,\dots,x_n)=1$. And we will also need
$s_j(x_1,\dots,x_n)=\sum_{i=1}^nx_i^j$.

\begin{lemma} $s_j(1,2,\dots,n-1)$ is a polynomial in $n$ with coefficients
in ${\Bbb Z}_{(j+2)}$, and $s_j(1,2,\dots,n-1)\equiv0$ mod $n$.
\end{lemma}

\begin{proof} By the Euler-Maclaurin formula for
$s_{j}$ [HW], we have
$$s_{j}(1,2,\dots , n-1)=\sum_{i=0}^{j} \frac{1}{j+1-i}\binom{j}{i}
n^{j+1-i}\beta_{i},$$
where $\beta_{0}=1, \beta_{1}=-1/2$, and
$\beta_{2i}=(-1)^{i-1}B_{i}$, $\beta_{2i+1}=0$ for
$i\geq 1$ ($B_{i}$ is the Bernoulli number).
We also have  von Staudt's theorem which determines the
fractional part of $B_{i}$ [HW].  It follows that $\beta_{i}\in
 {\Bbb Z}_{(i+2)}$.  This implies the lemma.
\end{proof}

For a pair of non-negative integers $l,m$, we define
\begin{equation}
H_{l,m}(r)=\frac{\hr!}{(\frac{r-1}{2}-l+m)!} (-1)^{m+1}
\sigma_{m}(1,2,\cdots , \frac{r-1}{2}-l+m-1).
\end{equation}

\begin{lemma} $H_{l,m}$ is a Fermat function. Moreover, let
$$g'_{l,m}=\r(H_{l,m}),$$
then $g'_{l,m}\in {\Bbb Z}_{(m+2)}[\frac12]$.
\end{lemma}

\noindent{\it Remark:} The proof of Lemma 2.4 will also provide us a
recursive formula to calculate $g'_{l,m}$. See Formula (2.5).

\begin{proof} We first derive a formula for $\sigma_{j}(1,2,\cdots ,
n-1)$ using a trick of Lagrange. We will abbreviate
$\sigma_j=\sigma_j(1,2,\dots,n-1)$.

We have $$
(x-1)(x-2)\cdots (x-(n-1))=x^{n-1}-\sigma_{1}x^{n-2} +\cdots
+(-1)^{n-1}\sigma_{n-1} .$$
If we multiply both sides by $x$ and change $x$ into $x-1$, we
obtain
$$\begin{aligned}
 (x-1)^{n}-&\sigma_{1}(x-1)^{n-1}+\cdots + (-1)^{n-1}\sigma_{n-1}(x-1)\\
&= (x-1)(x-2)\cdots (x-n)\\
&=(x-n) (x^{n-1}-\sigma_{1}x^{n-2}+\cdots +(-1)^{n-1}\sigma_{n-1}).
\end{aligned}
$$
Equating coefficients, we obtain
$$
\begin{aligned}
\sigma_{1}(1,2, \dots , n-1)&=\binom{n}{2}, \\
2\sigma_{2}(1,2,\dots,n-1)&=\binom{n}{3}+\binom{n-1}{2}\sigma_{1}, \\
3\sigma_{3}(1,2,\dots,n-1)&=\binom{n}{4}+\binom{n-1}{3}\sigma_{1} +
\binom{n-2}{2}\sigma_{2} ,\\
\dots &\dots \dots,\\
j\sigma_{j}(1,2,\dots,n-1)&=\sum_{i=0}^{j-1} \binom{n-i}{j+1-i}
\sigma_{i},\\
\dots &\dots \dots,\\
 (n-1) \sigma_{n-1}(1,2,\dots,n-1)&=
\sum_{i=0}^{n-2}\sigma_{i}.
\end{aligned}$$

Therefore, we see inductively that $\sigma_{j}(1,2,\dots, n-1)$ is
a polynomial of $n$ with coefficients in ${\Bbb Z}_{(j+2)}$ and divisible by
$n(n-1)\cdots (n-j)$.  Using these formulae for $\sigma_{j}$,
we have
\begin{equation}\begin{aligned}
&\frac{j\sigma_{j}(1,2,\dots ,n-1)}{n(n-1)\cdots (n-j)}\\
&=
\frac{1}{(j+1)!}+ \sum_{i=1}^{j-1}
\binom{n-i}{j+1-i}
\frac{\sigma_{i}(1,2,\dots ,n-1)}{n(n-1)\cdots (n-j)}\\
&=\frac{1}{(j+1)!}+ \sum_{i=1}^{j-1}
\frac{(n-i)}{(j+1-i)!}\frac{1}{i}
\frac{i \sigma_{i}(1,2,\dots ,n-1)}{n(n-1)\cdots (n-i)}.
\end{aligned}
\end{equation}

Suggested by (2.4), we define a series of rational
numbers $\sigma_{i}^{l,m}\in {\Bbb Z}_{(m+2)}[\frac12]$ for $i\geq 1$
recursively as follows:
$$\begin{cases}
\sigma_{1}^{l,m}=\frac{1}{2}, & \\
\sigma_{j}^{l,m}=\frac{1}{(j+1)!} +\sum_{i=1}^{j-1}
\frac{(-\frac12-l+m-i)}{(j+1-i)!}\cdot \frac{\sigma_{i}^{l,m}}{i}& \text{for $j\geq
2$.}
\end{cases}
$$
It is clear that we may identify $\sigma_{j}^{l,m}$ with
the residue of
$$\frac{j\sigma_{j}(1,2,\dots,n-1)}{n(n-1)\cdots (n-j)}$$
with
$$n=\frac{r-1}{2}-l+m.$$

Then we have
\begin{equation}
\begin{cases}
g_{0,0}'=-1& \\
g_{l,0}'=(-1)(-\frac12)\cdots (-\frac12-l+1)& \text{if $ l\geq 1$},\\
g_{l,m}'= (-1)^{m+1} (-\frac12)\cdots (-\frac12-l)
 \frac{\sigma_{m}^{l,m}}{m} & \text{if $m\geq 1$}.
\end{cases}
\end{equation}

Obviously, Formula (2.5) for $g_{0,0}'$ and $g_{l,0}'$ follows from
the formula of
$\text{Res}(D_k)$ in Example 2.1.
To see Formula (2.5) for
$m\geq 1$, we only need to notice that
$$\frac{(-\frac12)\cdots(-\frac12-l)}
{(-\frac12-l+m)\cdots(-\frac12-l)}=\text{Res}(D_{l-m}).$$
This finishes the proof of Lemma 2.4.
\end{proof}

\section{Fermat limits and Ohtsuki's invariants}

We fix some notation about the \lq\lq quantum parameter'' first:

\begin{itemize}

\item We will use $r$ to denote an odd
prime and $q=e^{2 \pi\sqrt{-1}/r}$.  The quantum integer
$$[k]=\frac{q^{\frac{k}2}-q^{-\frac{k}2}}{q^{\frac12}-q^{-\frac12}}.$$
Note that
$$[2]=q^{\frac12}+q^{-\frac12}.$$

\item For $R=\Bbb Z$ or ${\Bbb Z}_{(r)}$, $O((q-1)^k;R)$
stands for a complex number
of the form $u(q-1)^k$ for some $u\in R[q]$.

\item If $f(x)$ is a $C^{\infty}$
function, we use $\text{Coeff}_{n}(f;x-a)$ to denote
the coefficient of $(x-a)^{n}$ in the
Taylor expansion of $f$ at $x=a$.

\item The {\it Gauss sum} is $G_0(q)=\sum_{k=0}^{r-1}q^{k^2}$. The {\it
weighted Gauss sum} is $G_{2l}(q)=\sum_{k=0}^{r-1}k^{2l}q^{k^2}$.

\end{itemize}

To define Ohtsuki's invariants, we need the following definition.

Let $l(r)$ be an integer valued function of $r$ with
$$\lim_{r\rightarrow +\infty} l(r)=+\infty\qquad\text{and}
\qquad l(r)\leq r-2.$$
Given a sequence of complex numbers $s_r(q)\in {\Bbb Z}_{(r)}[q]$.
Fix a non-negative integer $n$; for
sufficiently large $r$, we write
$$
s_r(q)=a_{r,0}+a_{r,1}(q-1)+\cdots+a_{r,n}(q-1)^n+\cdots
            +a_{r,l(r)}(q-1)^{l(r)}+O((q-1)^{l(r)+1},{\Bbb Z}_{(r)})
$$
for some $a_{r,n}\in {\Bbb Z}_{(r)}$.

\begin{definition} The sequence of complex numbers $s_r(q)$ has
a {\it Fermat limit}, denoted by $\text{\em f-lim}\,s_r(q)$,
if each $a_{r,n}$, thought of as a function of $r$,
is a Fermat function. Let $\lambda_n=\r(a_{r,n})$. We write
$$\text{\em f-lim}\,s_r(q)=\sum_{n=0}^{\infty}\lambda_n(t-1)^n
\in {\Bbb Q}[[t-1]].$$
\end{definition}

\begin{lemma} The Fermat limit of $s_r(q)$ is well-defined and
if it exists, it is unique.
\end{lemma}

\begin{proof} Notice that although
the expansion $s_r(q)$ is not unique as $q$ is not a free variable,
$a_{r,n}\in{\Bbb Z}/r{\Bbb Z}$ is well defined for $0\leq n\leq l(r)$.
This is because of $r=O(\Q^{r-1};\Bbb Z)$, which follows from the
only relation
in ${\Bbb Z}_{(r)}[q]$: $\sum_{i=0}^{r-1} q^{i}=0$. Therefore, each residue
$\lambda_n$ is well-defined regardless of the choice of $a_{r,n}$.

The uniqueness of a Fermat limit now comes from the uniqueness of residues
of Fermat functions (Lemma 2.1).
\end{proof}

\noindent{\it Remark:} Notice that if $\lambda_n$ is an integer if and only if
the congruence class $\overline{a_{r,n}}\in{\Bbb Z}/r{\Bbb Z}$, which is
well-defined, stabilizes for $r$ sufficient large.

With all the above understood, we summarize the main results of [Mu] and [O1]
into the following theorem.

\begin{theorem}
Let $M$ be an oriented homology 3-sphere, and $\tau_{r}(M)$ be the
quantum $SU(2)$ invariant of Reshetikhin and Turaev {\em [RT]} at the
r-th root of unity normalized as in {\em [KM]}.  Then

\begin{enumerate}
\item {\em (H. Murakami)} $\tau_{r} (M) \in {\Bbb Z}[q]$;

\item {\em (T. Ohtsuki)} $\text{\em f-lim}\,\tau_{r} (M)=\sum_{n=0}^\infty
\lambda_{n} (t-1)^{n}$, and $\lambda_{n}\in {\Bbb Z}_{(2n+2)}$;

\item {\em (H. Murakami)} $\lambda_{0}=1$, and $\lambda_{1}=6
\lambda_{C}$, where $\lambda_{C}$ is the Casson invariant.
\end{enumerate}
\end{theorem}

It follows from Theorem 3.1 (2) and Lemma 3.1 that the $\lambda_{n}$
are topological
invariants of oriented integral homology 3-spheres, and we call
$\lambda_{n}$ the {\em n-th Ohtsuki invariant} of integral homology
3-spheres.

\section{Invariants of links}

Let $L=K_{1}\cup K_{2}\cup\cdots \cup K_{\# L}$ be an oriented link in
$S^{3}$, where we use $\# L$ to denote the number of components of $L$.
We will denote the unknot by $O$, and the empty link by
$\emptyset$.

The Conway polynomial $\nabla (L;z)\in{\Bbb Z}[z]$ is defined
by
$$\begin{cases} \nabla (O;z)=1 \\
\nabla(\emptyset;z)=0 \\
\nabla (L_{+};z)-\nabla (L_{-};z)=-z\cdot \nabla(L_{0};z).
\end{cases}
$$
Here, as usual, $L_{+}$, $L_{-}$ and $L_{0}$ are the links which have plane
projections identical to each other except
in one small disk where their projections are a positive crossing, a
negative crossing and an orientation preserving smoothing of that
crossing, respectively.
Note the negative sign on the right hand side of the skein relation.
Its effect is to change the usual Conway polynomial by a normalization factor.

The Jones polynomial $V(L;t)\in
{\Bbb Z}[t^{\frac12}, t^{-\frac12}]$ is defined by
$$\begin{cases}
V(O;t)=1;\\
V(\emptyset,t)=(t^{\frac12}+t^{-\frac12})^{-1}\\
tV(L_{+};t)-t^{-1}V(L_{-};t)=(t^{\frac12}-t^{-\frac12})V(L_{0};t).\end{cases}
$$
Note that our normalization differs
from the usual definition of the Jones polynomial in [J].  Actually it
is obtained from the usual Jones polynomial by changing $t$ to
$t^{-1}$ and multiplying by $(-1)^{\# L-1}$.  We put
$$X(L;t)=\frac{V(L;t)}{\T^{\# L -1}}.$$
Then
$$X(O;t)=X(\emptyset:t)=1.$$
We also put
$$\Phi(L;t)=\sum_{L'\subset
L}(-1)^{\# L-\# L'} X(L';t)$$
where the sum runs over all sublinks of $L$
including the empty link, and $L$ itself.
We have $\Phi(O;t)=0$ and define
$\Phi(\emptyset;t)=0$.

It will be useful to note that
$X(L;t),\Phi(L;t)\in{\Bbb Z}[t, t^{-1}, \frac{1}{t+1}]$.

The normalized Jones polynomial $X(L;t)$ and the \lq\lq averaged'' Jones
polynomial behave nicely under disjoint unions. Namely, if $L$
splits geometrically into two links $L_{1}$ and
$L_{2}$, then
$$X(L;t)=X(L_1;t)\cdot X(L_2;t)$$
and
$$\Phi(L;t)=\Phi(L_{1};t)\cdot \Phi(L_{2};t).$$

Further, we put
$$\Phi_{i}(L)=\left.\frac{d^{i}\Phi(L;t)}{dt^{i}}\right|_{t=1}$$
so that
$$\Phi(L;t)=\sum_{i=0}^\infty\frac{\Phi_i(L)}{i!}(t-1)^i.$$
Finally for each $i\geq 1$, we set
$$\phi_{i}(L)=\frac{(-2)^{\# L}}{(\# L+i)!}\cdot \Phi_{\# L +i}(L).$$
The invariants $\phi_i$ will be the basic link invariants we use
to express $\lambda_{n}$.

Given a link $L =K_{1}\cup K_{2}\cup \cdots \cup K_{\mu}$ with a certain projection.
Consider a pair of crossings between two components
$K_{i}$ and $K_{j}$, $i\neq j$,  with different signs. Denote this link by $L_{+-}$.
We may modify the link diagram at these two crossings to get various
other links, say $L_{-+}$,
$L_{--}$,  $L_{++}$,  $L_{0-}$, $L_{-0}$ and $L_{0,0}$.  We say
$L_{+-}$ is obtained from $L_{-+}$ by a {\em double crossing change.}
A link $L =K_{1}\cup K_{2}\cup \cdots \cup K_{\mu}$ is called an {\em algebraically
split link (ASL)} if the linking number between every pair of
components of $L$ is 0.  In particular, a knot is always an ASL.
If there are mutually disjoint 3-balls $B_{1}, B_{2},\dots ,
B_{m}$ in $S^{3}$ with $K_{i}\subset int B_{i}$, then $L$ is called
a {\em geometrically split link} (GSL).
Note that every algebraically split link can be
changed into a geometrically split link by a finite number of
double crossing changes, and the minimum number needed is called the
{\em double unlinking number} of this ASL.

\begin{lemma} {\em (Double crossing change formulae)}
For any link $L_{+-}$ with
$\# L_{+-}\geq 2$, we have
\begin{equation}
(t^{2}+t)[\Phi(L_{+-};t)-\Phi(L_{-+};t)]=(t-1)[
 \Phi(L_{-0};t)-\Phi(L_{0-};t)] .
\end{equation}
and
\begin{equation}
\begin{aligned}
\phi_{i} (L_{+-})- \phi_{i}(L_{-+})&=-[\phi_{i}(L_{0-}) -
\phi_{i}(L_{-0}) ]
\\
- \frac32\,[ \phi_{i-1} (L_{+-}) &-\phi_{i-1} (L_{-+}) ] -\frac12\,
[\phi_{i-2}(L_{+-}) -\phi_{i-2}(L_{-+}) ].
\end{aligned}
\end{equation}
\end{lemma}

\begin{proof}
Using the skein relation for $V(L;t)$, and
$\# L_{+-}=\#L_{-+}=\# L_{0-}+1=\#L_{-0}+1$, we have
$$(t^{2}+t)X(L_{+-};t)-(1+t^{-1})X(L_{--};t)=(t-1)X(L_{0-};t);$$ and
similarly,
$$(t^{2}+t)X(L_{-+};t)-(1+t^{-1})X(L_{--};t)=(t-1)X(L_{-0};t).$$
By subtracting the second equation from the first, we obtain
$$(t^{2}+t)[X(L_{+-};t)-X(L_{-+};t)]=(t-1)[
 X(L_{-0};t)-X(L_{0-};t)] .$$
By a direct computation, we get
$$(t^{2}+t)[\Phi(L_{+-};t)-\Phi(L_{-+};t)]=(t-1)[
 \Phi(L_{-0};t)-\Phi(L_{0-};t)] .$$
This is (4.1). Expanding both sides into Taylor series at $t=1$ and equating
coefficients, we will get (4.2) since $\# L_{+-}=\#L_{-+}=\# L_{0-}+1=\#L_{-0}+1$.
\end{proof}

The following lemma partially explains the definition of $\phi_i$. Notice that
Lemma 3.3 (2) is a special case of Proposition 3.4 in [O1] (see also Lemma 5.2),
proved here using the double crossing change formula.

\begin{lemma}  We have
\begin{enumerate}
\item $\frac{\Phi_{i}(L)}{i!} \in{\Bbb Z}[\frac12]$ for each $i$;
\item If $L$ is an ASL, then $\Phi_{i}(L)=0$ if $i\leq \#L$.
\end{enumerate}
\end{lemma}

\begin{proof}
(1) This follows from the fact that
$\Phi(L)\in{\Bbb Z}[t,t^{-1}, \frac{1}{t+1}]$.

(2) When $\#L=1$, (2) follows from the fact that
$$V(L,1)=1\qquad\text{and}\qquad\left.\frac{dV(L,t)}{dt}\right|_{t=1}=0.$$
The general case can be proved inductively using the double crossing
change formula (4.1) together with the fact that
$$\Phi(L;t)=\prod_{i=1}^{\#L}\Phi(K_i;t)$$
if $L=K_1\cup K_2\cup\cdots\cup K_{\#L}$ is a GSL.
\end{proof}

Let
$$v_{i}(L)=\left.\frac{d^{i}V(L;e^h)}{dh^{i}}\right|_{h=0}.$$

\begin{lemma}
Let $K$ be a knot.  Then
\begin{enumerate}

\item $v_{2}(K)  \in 6\Bbb Z.$

\item $v_{3}(K) \in 9\Bbb Z.$

\item $v_{4}(K) \in 6\Bbb Z.$

\item $v_{2}(K) + v_{4}(K)\in 18\Bbb Z.$
\end{enumerate}
\end{lemma}

\begin{proof} Let $K_+$ and $K_-$ be a pair of knots differed by a crossing
change.  As before, $K_0$ is the two component link obtained by an
orientation preserving smoothing of that crossing. We denote by $l$ the
linking number of the two components of $K_0$. We may have another knot
$K_{\infty}$ by smoothing the given crossing inconsistent with the
orientation. We have
\begin{equation}
V(K_+;t)-t^{-1}V(K_-;t)=(1-t^{-1})\,t^{-3l}\,V(K_\infty;t).
\end{equation}
See Corollary 13.4 in [J].

We use
$v_{i}$ with subscripts $+$, $-$ or
$\infty$ to denote the corresponding $v_{i}$ of $K_{+}$, $K_{-}$ and $K_{\infty}$.
By taking derivatives on both side of (4.3), we have:
$$\begin{aligned}
v_{2+}-v_{2-}&=-6l \\
v_{3+}-v_{3-}&=-3v_{2-}+27l^{2}+9l+3v_{2\infty}\\
v_{4+}-v_{4-}&=-4v_{3-}+6v_{2-}+4v_{3\infty}
-6v_{2\infty}(1+6l)-12l-54l-108l^{3}.\end{aligned}$$
Then the lemma follows inductively.
\end{proof}

\begin{theorem} Let $L$ be an ASL link.  Then
\begin{enumerate}
\item $\phi_{1}(L) \in 6\Bbb Z.$
\item $\phi_{2}(L) \in 3\Bbb Z .$
\end{enumerate}
\end{theorem}

\begin{proof}
(1) Using the double crossing change formula (4.2), we get
$$\phi_{1}(L_{+-}) - \phi_{1}(L_{-+})
=-[\phi_{1}(L_{0-}) - \phi_{1}(L_{-0})]=0.$$
By induction on the double unlinking number and the number of
components, it suffices to prove (1) for GSL's. If
$L$ is a knot, then
$\phi_{1}(L)=-v_{2}(L)\in6\Bbb Z$. If $L$ is a GSL with $\#L\leq2$, then
$\phi_2(L)=0$. So (1) holds for GSL's and, therefore, also for ASL's .

(2) Using the double crossing formula (4.2) again:
$$\phi_{2}(L_{+-}) - \phi_{2}(L_{-+})
=-[\phi_{2}(L_{0-}) - \phi_{2}(L_{-0})]
-\frac32\,[\phi_{1}(L_{+-}) - \phi_{1}(L_{-+})].$$
By induction on the double unlinking number and the number of
components, it suffices to prove (2) for GSL's.
If $L$ is a knot,
$$\phi_{2}(L)=-\frac{v_{3}(L)-3v_{2}(L)}{3}.$$
For a two-component GSL $L=K_{1}\cup K_{2}$,
$\phi_{2}(L)=v_{2}(K_{1})\cdot v_{2}(K_{2})$.  For a GSL $L$ with $\#L\geq3$,
$\phi_{2}(L)=0$. Thus (2) holds.
\end{proof}

In the light of Theorem 4.1, we was tempted to conjecture originally that
$n!\, \phi_{n}(L)\in6\Bbb Z$ for every ASL $L$. But this was shown by
Boden [Bo]
to be not true. We therefore strengthen the conjecture to the following form.

\begin{conjecture} For every boundary link $L$, $n!\, \phi_{n}(L)\in
6\Bbb Z.$
\end{conjecture}

\section{The second invariant $\lambda_{2}$}

A framed ASL is said to be {\em unit framed} if the framing of each component
is $\pm 1$.  Given a link $L$ and a positive integer $m$, we use
$L^{m}$ to denote the 0-framed $m$-parallel of $L$, i.e. each component in
$L$ is replaced by $m$ parallel copies having linking number zero with each
other.  So
$L^{m}$ is  an ASL if $L$ is.  Sublinks of $L^m$ will be assumed to be in one-one
correspondence with $\mu$-tuples $(i_1,\dots,i_\mu)$, where $\mu=\#L$,
in such a way that $L'$ has $i_\xi$ parallel copies of the $\xi$-th component of
$L$, $0\leq i_\xi\leq m$.
If $L$ is a framed link, $L^m$ and all its sublinks will inherit a framing from
$L$. If the $\xi$-th component of $L$ is framed by $f_\xi$,
$\xi=1,\dots,\mu$, we denote
$$f_L=\prod_{\xi=1}^{\mu}f_\xi.$$
As before, if $L$ is a framed link in $S^{3}$, then
$S^{3}_{L}$ denotes the resulting 3-manifold from the Dehn surgery on $L$.

\begin{theorem} Let $L$ be a unit framed ASL and $S^{3}_{L}$ be the
homology 3-sphere obtained from Dehn surgery on $L$. Let $\lambda_1$ and
$\lambda_2$ be the first and second Ohtsuki invariants, respectively. Then
\begin{equation}
\lambda_{1} (S^{3}_{L})=\sum_{L'\subset L}f_{L'}\,\phi_{1}(L')
\end{equation}
and
\begin{equation}
\lambda_{2} (S^{3}_{L})=\sum_{L'\subset L}\phi_{1}(L') \,
 f_{L'}\, \frac{ \# L'}{2}
+\sum_{L'\subset L^{2}}\phi_{2}(L')\,
f_{L'}\, \frac{1}{2^{s_{2}(L')}}.
\end{equation}
Here, if $L'$ corresponds to the $\mu$-tuple $(i_1,\dots,i_\mu)$, $s_{2}(L')=
\#\{i_\xi\,;\,i_\xi=2\}$.
\end{theorem}

The formula for $\lambda_{1}$ is equivalent to Hoste's formula
for the Casson invariant [Ho].  This can be proved using Corollary 3.12 [Mu].
Both formulae will be proved in Section 6.

Combining with Theorem 4.1, we obtain:

\begin{corollary}
Let $M$ be an oriented homology 3-sphere.  Then
\begin{enumerate}
\item $\lambda_{1}(M) \in 6\Bbb Z.$
\item $\lambda_{2}(M) \in 3{\Bbb Z}[\frac12] .$
\end{enumerate}
\end{corollary}

Notice that Corollary 4.1 (2) will be strengthened in Theorem 7.1.

In the rest of this section, we will study the behaviors of the invariant
$\lambda_2$ on homology 3-spheres obtained from $1/n$-Dehn surgery on
knots. We will use (5.2) to make the relation between $\lambda_2$ with
(coefficients of) the Jones polynomial very explicit in this particular case.

\begin{theorem}
Let $n$ be an integer, and
$S^{3}_{K,1/n}$ be the homology 3-sphere obtained from $1/n$-Dehn
surgery on a knot $K$.  Then
\begin{equation}
 \lambda_{2}(S^{3}_{K,1/n})=\frac{n}{2}\, v_{2}(K)
-\frac{n}{3}\, v_{3}(K) -\frac{n^{2}}{6}\, [v_{2}(K) +v_{4}(K) -4 v_{2}^{2}(K)].
\end{equation}
\end{theorem}

To prove this theorem, we need two lemmas.

\begin{lemma}
\begin{enumerate}
\item $\Phi_{2}(K)=v_{2}(K)$, and $\Phi_{3}(K)=v_{3}(K)-3v_{2}(K).$
\item $\Phi_{4}(K^{2})=-2v_{2}(K)-2v_{4}(K) + 8 v_{2}^{2}(K).$
\end{enumerate}
\end{lemma}

\begin{proof} (1) Exercise.

(2) $\Phi_{4}(K^{2})$ is a Vassiliev invariant
of order 4 (see, for example, [MM]).  Therefore, for any knot $K$,
$$\Phi_{4}(K^{2})=a\cdot v_{2}(K) + b\cdot v_{3}(K) + c\cdot
v_{4}(K) + d \cdot \phi(K) + e \cdot v_{2}^{2}(K),$$ where
$\phi$ is the coefficient of $z^{4}$ in the Conway polynomial
of $K$, and $a,b,c,d,e$ are some constants.  If $K$ is a torus knot,
then the Conway/Jones polynomials of $K$ and $K^{2}$ are known
(see, for example, [M]).   Using these data, we can determine
$a,b,c,d,e$.  And this gives the above formula.
\end{proof}

\begin{lemma}
Let $L$ be an ASL with $\mu=\#L$. If $L'\subset L^m$ corresponding to
$(i_1,\dots,i_\mu)$ with a certain
$i_\xi=m$, then $\Phi_i(L')=0$ for $i\leq\mu+m-1$.
\end{lemma}

Compare it with Lemma 4.2 (2). This is the result of Ohtsuki mentioned
above. It is proved in [O1] using the colored Jones polynomial.

\begin{proof-5-2} Assume $n\neq0$. Using Kirby's calculus on framed links,
one can see that
$S^{3}_{K,1/n}$ is the same as $S^{3}_{K^{|n|}}$ with the framing
of each component equal to $\text{sign}(n)$. This is done explicitly in [G].  By
(5.2) and Lemma 5.2, we have
$$\lambda_{2}(S^{3}_{K_,1/n})=-\frac{n}{2} \cdot \Phi_{2}(K) -\frac{n}{3}
\cdot \Phi_{3}(K) +\frac{n^{2}}{12}\cdot \Phi_{4}(K^{2}).$$
Then (5.3) follows from Lemma 5.1.
\end{proof-5-2}

Combining (5.3) with Lemma 4.4, we get

\begin{corollary}
Let $M$ be the resulting homology 3-sphere
of 1/n-Dehn surgery on a knot $K$.  Then
 $\lambda_{2}$ is an integer divisible by 3, i.e., $\lambda_{2}\in 3\Bbb Z$.
\end{corollary}

Once again, this conclusion will be strengthened in
Theorem 7.1.

\begin{corollary}
If $v_{3}(K)\neq 0$ for a knot K, and let $K^{*}$ be the mirror
image of $K$.  Then $S^{3}_{K,1/n}$ and $S^{3}_{K^{*},1/n}$
are distinct as unoriented manifolds for each $n\neq 0$.
\end{corollary}

\begin{proof}
It suffices to prove that $S^{3}_{K,1/n}$ is neither homeomorphic to
 $S^{3}_{K^{*},1/n}$ nor
 $\overline{S^{3}_{K^{*},1/n}}$.

If $v_2(K)\neq 0$, then $\lambda_1$ will show this.  If $v_2(K)=0$,
then $\lambda_2$ takes the same values on
 $S^{3}_{K^{*},1/n}$ and
 $\overline{S^{3}_{K^{*},1/n}}$ (see the remark below).
As $v_{2n}(K)=v_{2n}(K^{*})$,
and  $v_{2n+1}(K)=-v_{2n+1}(K^{*})$, so
$$\lambda_{2}(S^{3}_{K,1/n})
-\lambda_{2}(S^{3}_{K^{*},1/n}) =-\frac{2n}{3}\, v_{3}(K).$$
This completes the proof.
\end{proof}

\begin{corollary}
Let $n,m$ be two distinct integers, and $K$ be a knot.
\begin{enumerate}
\item If $v_{2}(K)\neq 0$,
then $S^{3}_{K,1/n}$ and
$S^{3}_{K,1/m}$ are distinct as unoriented 3-manifolds.
\item If $v_{2}(K)=0$, and $v_{3}(K)\neq
-\frac{n+m}{2} \cdot
 v_{4}(K) $, then
 $S^{3}_{K,1/n}$ and
$S^{3}_{K,1/m}$ are distinct as unoriented 3-manifolds.
\end{enumerate}
\end{corollary}

\begin{proof}
(1) follows from comparing $\lambda_{1}$ of the two manifolds.

(2) follows from comparing $\lambda_{2}$ of the two manifolds.
\end{proof}

Some examples are as follows.

\begin{example}{\em
(1) If $M$ is the Poincare homology 3-sphere $\Sigma(2,3,5)$ ($+1$-surgery
on the right-handed trefoil knot), then $\lambda_{2}(M)=39$.

(2) If $M$ is the homology 3-sphere
$\Sigma(2,3,7)$ ($+1$-surgery on the
left-handed trefoil knot), then $\lambda_{2}(M)=63$.}
\end{example}

\noindent{\it Remark:}  Suppose $M$ is an oriented
homology 3-sphere and $\overline{M}$ is
the oriented homology 3-sphere obtained from $M$ by reversing the
orientation.  Then $\lambda_{2} (\overline{M})=\lambda_{2}(M)+ \lambda_{1}(M).$
If $N$ is the oriented homology 3-sphere obtained from (-1)-surgery
on a knot $K$ and $M$ is the oriented homology 3-sphere obtained from
(+1)-surgery on the mirror image of $K$, then $N=\overline{M}$.  Using
these facts, we observe that our computation above agrees completely
with that in [L].

\begin{example}{\em
If $M$ is the homology 3-sphere
$\Sigma(2,5,7)$ (this is the boundary of the
{\it Mazur manifold}), then
$\lambda_{2}(M)=-66$. Notice that the Mazur manifold is contractible and
$6|\lambda_2$ in this particular case. This raises the question of whether
$\lambda_2/3$ mod 2 is a homology spin cobordism invariant.}
\end{example}

It is not very difficult to prove that there exists
two knots $K_{1}, K_{2}$ with the same Jones polynomial and
such that $1/n$-Dehn surgery on $K_{1}$ and $K_{2}$, respectively,
for some $n$ give distinct homology 3-spheres.
If $M$ is obtained from Dehn surgery on a knot,
it follows from Theorem 5.2 that $\lambda_{2}(M)$ is determined by the Jones
polynomial. So $\lambda_{2}(M)$ can be the same for distinct
homology 3-spheres.

\section{The Formula for $\lambda_{n}$}

Let $M$ be the resulting oriented homology 3-sphere of Dehn
surgery on a unit framed ASL $L$ and $\mu=\#L$.  Let $f_\xi=\pm1$ be the
framing on the $\xi$-th component of $L$.  Then for each
$n\geq 1$, Ohtsuki gave the following expression for $\lambda_{n}$:
$$
\begin{aligned}
&\lambda_n(M):=\\
&\text{Coeff}_{n}\left\{q^{\frac{3}{4}
\sum f_{\xi}-\frac{\mu}{2}}\sum_{l=1}^{n}\sum_{L'\subset L^l}
\frac{\Phi_{(l+\# L')}(L')}{(l+\# L')!} \Q^{l} \prod_{\xi=1}^{\mu}
\left(-f_{\xi} \sum_{m_{\xi}=0}^{n-l}h_{f_\xi,i_{\xi},m_{\xi}}{\Q}^{m_{\xi}}\right);
q-1 \right\},
\end{aligned}$$  where $h_{f_\xi, i_{\xi}, m_{\xi}}$'s are some unknown
constants whose existence was established by Ohtsuki.
Our main result here is
a more practical formula for $\lambda_{n}$.
In the following, we will restate several lemmas in [O1] using
Definition 3.1.  This unifies those lemmas and  clarifies the uniqueness of
the constants appearing there.

Let the constants $\nu_{f,i,m}$ be defined by the following formal
equation, i.e. by equating coefficients on both sides in the
expansion of power series in $(t-1)$:
$$ \sum_{m=0}^{+\infty} \nu_{f,i,m}
(t-1)^{m}=t^{\frac{3}{4}f-\frac{1}{2}}\cdot \sum_{m=0}^{+\infty}h_{f,i,m}
(t-1)^{m}.$$

\begin{theorem}
Let $M$ be as above.  Then
 $$\lambda_{n}(M)=\sum_{l=1}^{n}\sum_{L'\subset L^{l}}
\frac{\phi_{l}(L')}{(-2)^{\#L'}}\left\{ \prod_{\xi=1}^{\mu}(-f_{\xi})
\left(\sum_{m_{1}+m_{2}+ \cdots +m_{\mu}=n-l} \prod_{\xi=1}^{\mu}
\nu_{f_\xi,i_{\xi},m_{\xi}}\right)\right\}.$$
\end{theorem}

\begin{proof}
Using Ohtsuki's expression, we have
$$
\begin{aligned}
&\lambda_{n}(M)\\
&=\text{Coeff}_{n}\left\{\sum_{l=1}^{n}\sum_{L'\subset L^l}
\frac{\phi_l(L')}{(-2)^{\# L'}} \Q^{l} \prod_{\xi=1}^{\mu}
\left(-f_{\xi} q^{\frac{3}{4} f_\xi-\frac{1}{2}} \sum_{m_{\xi}=0}^{n-l}
h_{f_\xi,i_{\xi},m_{\xi}}{\Q}^{m_{\xi}}\right); q-1 \right\}\\
&=\text{Coeff}_{n}\left\{\sum_{l=1}^{n}\sum_{L'\subset L^l}
\frac{\phi_l(L')}{(-2)^{\# L'}} \Q^{l} \prod_{\xi=1}^{\mu}
\left(-f_{\xi} q^{\frac{3}{4} f_\xi-\frac{1}{2}}
 \sum_{m_{\xi}=0}^{+\infty}h_{f_\xi,i_{\xi},m_{\xi}}{\Q}^{m_{\xi}}\right);
q-1 \right\}\\
&=\sum_{l=1}^{n}\sum_{L'\subset L^l}
\frac{\phi_l(L')}{(-2)^{\# L'}} \text{Coeff}_{n-l} \left\{\prod_{\xi=1}^{\mu}
\left(-f_{\xi} q^{\frac{3}{4} f_\xi-\frac{1}{2}}
 \sum_{m_{\xi}=0}^{+\infty}h_{f_\xi,i_{\xi},m_{\xi}}{\Q}^{m_{\xi}}\right);
q-1 \right\}\\
&=\sum_{l=1}^{n}\sum_{L'\subset L^l}
\frac{\phi_l(L')}{(-2)^{\# L'}} \text{Coeff}_{n-l}\left\{\prod_{\xi=1}^{\mu}
\left(-f_{\xi}
 \sum_{m_{\xi}=0}^{+\infty}\nu_{f_\xi,i_{\xi},m_{\xi}}{\Q}^{m_{\xi}}\right);
q-1 \right\}.
\end{aligned}
$$
By picking out the $n$-th coefficient, we get the desired
formula.
\end{proof}

\begin{corollary} Let $K$ be a knot.  The knot invariant $\theta_n$
induced by $\lambda_{n}$ is defined to be
$\theta_n(K)=\lambda_{n}(S^3_{K,1})$.  Then $\theta_n$
is a Vassiliev invariant of order $2n$.
\end{corollary}

\begin{proof}
By Theorem 6.1, $\theta_n$ is a linear combination of
$\phi_{l}(K^{j})$ for $1\leq l\leq n, 1\leq j \leq l.$
As $\phi_{l}(K^{j})$ is a Vassiliev invariant of order $l+j$, $\theta_n$ is
an invariant of order $2n$.
\end{proof}

If $\lambda_{n}$ is indeed a finite invariant of order $3n$ in the
sense of [O2] as suggested by [R1,2], then Corollary 6.2 will follow from
the result in [Ha].

To get an explicit formula for $\lambda_{n}$, we have to determine
the constants $\nu_{f,i,m}$.  Our discussion here is parallel to
the discussion of $h_{f,i,m}$ in Section 8 of [O1].  To avoid
repetition of many formulae, we will quote some formulae from [O1]
and leave the
interested reader to consult [O1] for details. We will make some improvement
to his main
lemma (Lemma 8.3 in [O1]) and the result is more explicit.

To get $\nu_{f,i,m}$, fix a nonnegative integer $i$, and
let $l(r)= \frac{r-1}{2}-i-1$. Let
$$
\begin{aligned}
s_{r, f, i}(q)&=\left(\frac{f}{r}\right)
\cdot q^{\frac{3}{4} f-\frac{1}{2}}\cdot
\frac{\Q^{i+1}}{G_{0}(q)}\\
&\times \left\{ \sum_{k=1}^{\frac{r-1}{2}}
q^{\frac14f(k^{2}-1)} [k]\sum_{j=0}^{\frac{k-1}{2}} (-1)^{j}
\binom{k-j-1}{j}\binom{k-2j-1}{i}
[2]^{k-2j-1} \right\},
\end{aligned}$$
where $G_{0}(q)$ is the Gauss sum. Note that
$s_{r,f, i}\in {\Bbb Z}_{(r)}[q]$ and $\left(\frac{f}{r}\right)$ is the Legendre
symbol.

\begin{theorem} With respect to $l(r)$,
$$\text{\em f-lim}\, s_{r,f,i}(q)=\sum_{m=0}^{+\infty}
\nu_{f,i,m}(t-1)^{m}.$$
\end{theorem}

This is essentially Proposition 3.6 in [O1]. We need to compute the
Fermat limit of $s_{r,f,i}(q)$ in order to determine $\nu_{f,i,m}$.

Recall that $g'_{l,m}\in
{\Bbb Z}_{(m+2)}[\frac12]$ is the residue of $H_{l,m}$ given by (2.5).
Let $G_{2l}(q)$ be the
weighted Gauss sum.

\begin{lemma}
With respect to $l(r)=\frac{r-3}{2}$,  we have
$$\text{\em f-lim}\, \frac{\left(\frac{r-1}{2}\right)!\cdot G_{2l}(q)}
{(q-1)^{\frac{r-1}{2}-l}}
=\sum_{m=0}^{+\infty} g'_{l,m}(t-1)^{m}.$$
\end{lemma}

This is the main technical lemma and
is an improvement to Lemma 8.3 in [O1].
In [O1], Ohtsuki established the existence of $g'_{l,m}$, but his
proof does not give a formula for these numbers.  Here, the $g'_{l,m}$ are given
explicitly in (2.5).  We also have
$g'_{l,m} \in {\Bbb Z}_{(m+2)}[\frac12]$, whereas in [O1], it is only known that
$g'_{l,m} \in {\Bbb Z}_{(2m+2)}[\frac12]$.

\begin{proof}
To expand $G_{2l}(q)$ into a power series in $\Q$, we take the Taylor
expansion of $\sum_{k=0}^{r-1} k^{2l}x^{k^{2}}$ at $x=1$ and put
$x=q$. Hence we have $$
\begin{aligned}
&\left(\frac{r-1}{2}\right)!\, G_{2l}(q)\\
& =\sum_{m=1}^{r-2-l} \left(\frac{r-1}{2}\right)! \sum_{k=0}^{r-1}
\frac{k^{2l} k^{2} (k^{2}-1)\cdots (k^{2}-(m-1))}{m!} \Q^{m} \\
&\qquad\qquad\qquad +O(\Q^{r-1-l};\Bbb Z)\\
& =\sum_{m=1}^{r-2-l} \frac{\hr!}{m!} \sum_{k=0}^{r-1}
\sum_{1\leq m'\leq m}(-1)^{m-m'}\sigma_{m-m'}(1,\dots,m-1)k^{2l+2m'}
\Q^{m} \\
&\qquad\qquad\qquad +O(\Q^{r-1-l};\Bbb Z) \\
& =\sum_{m=1}^{r-2-l} \frac{\hr!}{m!}
\sum_{1\leq m'\leq m} (-1)^{m-m'}  \sigma_{m-m'}(1,\dots,m-1)
s_{2l+2m'}(1,\dots,r-1)  \Q^{m} \\
&\qquad\qquad\qquad +O(\Q^{r-1-l};\Bbb Z)
\end{aligned}
$$

If $2l+2m' \leq r-2$, i.e.,  $m'\leq m\leq \frac{r-3}{2}-l$, then
it follows from Lemma 2.3 and $r=O(\Q^{r-1};\Bbb Z)$ that
$$s_{2l+2m'}(1,\dots,r-1)=O(\Q^{r-1}; {\Bbb Z}_{(r)}).$$
Therefore,
$$\begin{aligned}
&\left(\frac{r-1}{2}\right)!\,G_{2l}(q)\\
&= \sum_{m=\frac{r-1}{2}-l}^{r-2-l} \frac{\hr!}{m!}
\sum_{1\leq m'\leq m} (-1)^{m-m'}  \sigma_{m-m'}(1,\dots,m-1)
s_{2l+2m'}(1,\dots,r-1)  \Q^{m}\\
&\qquad\qquad\qquad+O(\Q^{r-1-l};{\Bbb Z}_{(r)})
\end{aligned}
$$

We will use Fermat's Little Theorem ($a^{r-1}\equiv1$ mod $r$) to get rid of
$s_{2l+2m'}(1,\dots,r-1)$ in the above expression.
So if $2l+2m'\neq 0\; mod\; (r-1)$, then we may assume
$2l+2m'\leq r-2$. We have dropped these terms by the preceding argument.
So the only non-trivial contribution of $s_{2l+2m'}(1,\dots,r-1)$
comes from $2l+2m'\equiv0$ mod
$(r-1)$.  Since $\frac{r-1}{2}-l \leq m\leq r-2-l$, we have
$2l+2m'=r-1$, i.e. $m'=\frac{r-1}{2}-l$.  In this case,
$$s_{2l+2m'}(1,\dots,r-1)=\sum_{k=1}^{r-1} k^{r-1}\equiv-1\,\,\,\text{mod}\,\,r.$$
Hence we get
$$
\begin{aligned}
&\left(\frac{r-1}{2}\right)!\, G_{2l}(q)\\
&= \sum_{m=\frac{r-1}{2}-l}^{r-2-l} \frac{\hr!}{m!}
 (-1)^{m-\frac{r-1}{2}-l+1}
 \sigma_{m-\frac{r-1}{2}-l}(1,\dots,m-1) \Q^{m} \\
&\qquad\qquad\qquad +O(\Q^{r-1-l};{\Bbb Z}_{(r)}) \\
&= \sum_{m=0}^{\frac{r-3}{2}} \frac{\hr!}{(\frac{r-1}{2}-l+m)!}
(-1)^{m+1} \sigma_{m}(1,\dots,\frac{r-1}{2}-l+m-1)
 \Q^{\frac{r-1}{2}-l+m}\\
&\qquad\qquad\qquad+O(\Q^{r-1-l};{\Bbb Z}_{(r)})\\
&=\sum_{m=0}^{\frac{r-3}{2}} H_{l,m}(r)(q-1)^{\frac{r-1}{2}-l+m}.
\end{aligned}$$
This finishes the proof.
\end{proof}

For a pair of non-negative integers $l,m$,
let $g_{l,m}$ be a sequence of rational
numbers defined by the equation
$$\sum_{m=0}^{+\infty}g_{l,m}'(t-1)^{m}=\sum_{m=0}^{+\infty}g_{l,m}(t-1)^{m}
\sum_{m=0}^{+\infty}g_{0,m}'(t-1)^{m}.$$
Comparing the coefficients, we get
\begin{equation}
\begin{cases}
g_{l,0}= -g'_{l,0},& \\
g_{l,m}=-g'_{l,m}+\sum_{k=0}^{m-1} g_{l,k}g'_{0,m-k}& \text{if $m>0$.}
\end{cases}
\end{equation}
So we can compute $g_{l,m}$ recursively once the $g'_{l,m}$ are known.
And we also see that
$$g_{l,m}\in{\Bbb Z}_{(m+2)}[\frac12].$$

The following lemma is now a direct consequence of the definition.

\begin{lemma}
With respect to $l(r)=\frac{r-3}{2}$, we have
$$\text{\em f-lim}\, \Q^{l}\,\frac{ G_{2l}(q)}{G_{0}(q)}
=\sum_{m=0}^{+\infty}g_{l,m}(t-1)^{m}.$$
In particular, $g_{0,m}=\delta_{0,m}$.
\end{lemma}

Here comes our last technical lemma, whose proof is straightforward.

\begin{lemma}
Let $F_{i,l}(x)$ be the polynomial in $x$ defined recursively for
$i\geq 0, l\geq -1$ and $i\geq l$ by:
$$ F_{0,0}(x)=1,\,\, F_{i, -1}(x)=F_{i, i+1}=0,$$
$$ F_{i+1, l}(x)=F_{i, l-1}(x)-(2i+1-l)\, x\,F_{i,l}(x)+
(x^{2}-4)F'_{i,l}(x).$$  Then
\begin{enumerate}
\item $ F_{i,i}(x)=1$;

\item $F_{i+k, i}(x)$ is of the same degree $k$ for any $i\geq 0$; and

\item $F_{i+1,i}(x)=-x\,\binom{i+2}{2}.$
\end{enumerate}
\end{lemma}

Now we come to the computation of the Fermat limit of $s_{r,f,i}(q)$,
or equivalently, the determination of $\nu_{f,i,m}$.
First, note that $s_{r,f,i}(q)$ is the left-hand side of the equation in
Proposition 3.6 of [O1] multiplied by $\left(\frac{f}{r}\right)
\, q^{\frac34 f-\frac12}$. Then the
expansion of $s_{r,f,i}(q)$ into
 a power series in $\Q$
follows from the formulae of Ohtsuki on pages 108-109 of [O1].
There is an omission of $f^{l}$ on lines 12, 15 and 18 in his formulae
on page 109.
For $i=0$, we have $s_{r,f,i}=-f$.  So we assume $i\geq1$.  Then
we have
$$
\begin{aligned}
s_{r,f,i}(q)&=\left.
\frac{[2]^{i}q^{i+\frac12+\frac{f}{2}}} {i!}\right\{
\frac{F_{i,0}([2])}{\Q^{i+1}}
(q^{-f}-1) \\
& + q^{-f}\sum_{l=1}^{i}f^{l}F_{i,l}([2])(-2)^{l}q^{-\frac{l}2}
\sum_{l'=0}^{\frac{l}2}
\binom{l}{2l'}
f^{l'}\sum_{m=0}^{\frac{r-3}{2}}g_{l',m}\Q^{m-l'-i-1+l}\\
&\left.-\sum_{l=1}^{\frac{i}2}F_{i,2l}([2]) 2^{2l}q^{-l}f^{l}
\sum_{m=0}^{\frac{r-3}{2}}g_{l,m}
\Q^{m+l-i-1} \right\} +O(\Q^{\frac{r-1}{2}-i}; {\Bbb Z}_{(r)}).
\end{aligned}
$$
Note that $q^{-f}-1=(-f)\, q^{-\frac{f+1}{2}}(q-1)$ and
$[2]^{i}=q^{-\frac{i}2}(q+1)^{i}$. Separating $l'=0$ from the
second term,
we get
$$
\begin{aligned}
&s_{r,f,i}(q)= \frac{(-f)}{i!}\, q^{\frac{i}2}\,
(q+1)^{i}\,F_{i,0}([2]) (q-1)^{-i} \\
& + \sum_{l=1}^{i}f^{l}\,
\frac{(-2)^{l}}{i!}\, q^{\frac{i}2+\frac12-\frac{f}2-\frac{l}2}\,
(q+1)^{i}\,
F_{i,l}([2]) (q-1)^{-i-1+l}\\
& + \sum_{l=1}^{i}f^{l}\,
\frac{(-2)^{l}}{i!}\, q^{\frac{i}2+\frac12-\frac{f}2-\frac{l}2}\,
 (q+1)^{i}\,
F_{i,l}([2])
 \sum_{l'=1}^{\frac{l}2}\binom{l}{2l'}
f^{l'}\sum_{m=0}^{\frac{r-3}{2}}g_{l',m}\Q^{m-l'-i-1+l}\\
&-\sum_{l=1}^{\frac{i}2} f^{l}\,
\frac{2^{2l}}{i!}\, q^{\frac{i}2+\frac12+\frac{f}2-l}\, (q+1)^{i}
F_{i,2l}([2])
\sum_{m=0}^{\frac{r-3}{2}}g_{l,m}
\Q^{m+l-i-1} +O(\Q^{\frac{r-1}{2}-i}; {\Bbb Z}_{(r)}).
\end{aligned}
$$
By changing the order of summations in the third term, and combining
the sum with $l=2l'$ in the new order of summations with the fourth term,
we get
$$
\begin{aligned}
&s_{r,f,i}(q)=\frac{(-f)}{i!}\, q^{\frac{i}2}\,
(q+1)^{i}\, F_{i,0}([2]) (q-1)^{-i} \\
& + \sum_{l=1}^{i}f^{l}\,\frac{(-2)^{l}}{i!}\,
q^{\frac{i}2+\frac12-\frac{f}2-\frac{l}2}\, (q+1)^{i}\,
F_{i,l}([2]) (q-1)^{-i-1+l}\\
& + \sum_{l'=1}^{\frac{i}2}\sum_{l=2l'+1}^{i}f^{l+l'}\,
\frac{(-2)^{l}}{i!}\,\binom{l}{2l'}
\, q^{\frac{i}2+\frac12-\frac{f}2-\frac{l}2}\, (q+1)^{i}\,
F_{i,l}([2])
\sum_{m=0}^{\frac{r-3}{2}}g_{l',m}\Q^{m-l'-i-1+l}\\
&+ (q^{-f}-1)\, \sum_{l=1}^{\frac{i}2} f^{l}\,
\frac{2^{2l}}{i!}\, q^{\frac{i}2+\frac12+\frac{f}2-l}\, (q+1)^{i}
F_{i,2l}([2])
\sum_{m=0}^{\frac{r-3}{2}}g_{l,m}
\Q^{m+l-i-1}\\
&\qquad\qquad\qquad +O(\Q^{\frac{r-1}{2}-i}; {\Bbb Z}_{(r)})
\end{aligned}$$
$$
\begin{aligned}
&=\frac{(-f)}{i!}\, q^{\frac{i}2}\,
(q+1)^{i}\,  F_{i,0}([2]) (q-1)^{-i} \\
& + \sum_{l=1}^{i}f^{l}\,
\frac{(-2)^{l}}{i!}\, q^{\frac{i}2+\frac12-\frac{f}2-\frac{l}2}\,
 (q+1)^{i}\,
F_{i,l}([2]) (q-1)^{-i-1+l}\\
& + \sum_{l'=1}^{\frac{i}2}\sum_{l=2l'+1}^{i}f^{l+l'}\,
\frac{(-2)^{l}}{i!}\,\binom{l}{2l'}
\, q^{\frac{i}2+\frac12-\frac{f}2-\frac{l}2}\, (q+1)^{i}\,
F_{i,l}([2])
\sum_{m=0}^{\frac{r-3}{2}}g_{l',m}\Q^{m-l'-i-1+l}\\
&+ \sum_{l=1}^{\frac{i}2} (-f)\, f^{l}\,
\frac{2^{2l}}{i!}\, q^{\frac{i}2-l}\, (q+1)^{i}
\, F_{i,2l}([2])\,
\sum_{m=0}^{\frac{r-3}{2}}g_{l,m}
\Q^{m+l-i}\\
&\qquad\qquad\qquad +O(\Q^{\frac{r-1}{2}-i}; {\Bbb Z}_{(r)}).
\end{aligned}
$$
Using the fact that $g_{0,m}=\delta_{0m}$, we see that
the first term is the case $l=0$ of the fourth, and the second
is $l'=0$ of the third.
Combining these four terms into two terms,
and picking out the coefficient of $\Q^{m}$ in the expansion above,
we obtain
\begin{equation}
\begin{aligned}
&\nu_{f,i,m}=\sum_{l'=0}^{\frac{i}2}\sum_{l=2l'+1}^{i}f^{l+l'}\,
\frac{(-2)^{l}}{i!}\,\binom{l}{2l'}\,\text{Coeff}_{m} (t^{\frac{i+1-f-l}{2}}\, (t+1)^{i}\,F_{i,l}\T
\times\\
&\qquad\qquad \sum_{m'=0}^{+\infty}g_{l',m'}(t-1)^{m'-l'-i-1+l} ; t-1) \\
&\qquad+ \sum_{l'=0}^{\frac{i}2} (-f)\, f^{l'}\,
\frac{2^{2l'}}{i!}\,\text{Coeff}_{m} ( t^{\frac{i}2-l'}\, (t+1)^{i}
\, F_{i,2l'}\T \times\\
&\qquad\qquad\sum_{m'=0}^{+\infty}g_{l',m'}
(t-1)^{m'+l'-i}; t-1).
\end{aligned}
\end{equation}
As $g_{l,m}$ and $F_{i,l}(x)$ both can be computed recursively,
 so we can compute $\nu_{f,i,m}$.

\begin{theorem}
Let $M$ be an oriented integral homology 3-sphere. Then
$ \lambda_{n}(M) \in {\Bbb Z}_{(n+1)}$.
\end{theorem}
\begin{proof}
By Theorem 6.1,  we need $\nu_{f,i,m}$ with $0\leq m\leq n-1,$
and $0\leq i\leq n-m$ for $\lambda_{n}$.  Using (6.2),
we need $m'-l'-i-1+l\leq
m$ and
$m'+l'-i\leq m$. In both cases, the terms with $l'=0$
have the right denominator.  So we may assume $l'\geq 1$.
Then it is easy to check that $m'\leq n-1$.  So for
the formula of $\lambda_{n}$, we need $g_{l',m'}$
with $l'\leq n/2,$ and $m'\leq n-1$.
As $g_{l',m'}\in {\Bbb Z}_{(m'+2)}$, our theorem follows.
\end{proof}

This theorem tells us that the biggest factor in the denominator
of $\lambda_{n}$ is indeed $n$.  So it agrees with Conjecture 1.1.

To tie up the loose ends in Section 5, we prove the formulae for
$\lambda_{1}$ and $\lambda_{2}$.

\begin{proof-5-1}
First we need the following results for $\nu_{f,i,m}$ for $i=0,1,2$. Use (6.2),
for $i=0$
$$\nu_{f,0,m}=\text{Coeff}_{m}(-f; t-1),$$
so
$\nu_{f,0,0}=-f$, $\nu_{f,0,m}=0$ if $m\geq 1$. For $i=1$,
$$\nu_{f,1,m}=\text{Coeff}_{m}(t+1; t-1),$$
so $\nu_{f,1,0}=2$, $\nu_{f,1,1}=1$, and $\nu_{f,1,m}=0$ if $m\geq 2$. For
$i=2$,
$$\nu_{f,2,m}=\text{Coeff}_{m}\left(-\frac{t^{2}+2t+2}{2} (1+f +
4\sum_{m=1}^{+\infty}g_{1,m}(t-1)^{m-1}); t-1\right).$$
It follows that $\nu_{f,2,0}=-(2+2f +8 g_{1,1})$.  As
$g_{1,1}=-\frac14$,
so $\nu_{f,2,0}=-2f$.

For the formula for $\lambda_{1}$, we have $n=1$ and $l=1$.
Therefore,
$m_{1}=\cdots =m_{\mu}=0$
in Theorem 6.1.
Using the values of $\nu_{f,0,0}$, and
$\nu_{f,1,0}$, we get
$$\sum_{m_{1}+\cdots +m_{\mu}=0}
 \prod_{\xi=1}^{\mu}\nu_{f_\xi,i_{\xi},m_{\xi}}
=\prod_{\xi:i_{\xi}=0}
(-f_{\xi}) \, \prod_{\xi:i_{\xi}=1} 2.$$
This verifies the formula for $\lambda_1$.

For the formula of $\lambda_{2}$, we have
two cases: $n=2,\, l=1$ and $n=2,\, l=2$.
If $n=2,\, l=1$, then
$$\sum_{m_{1}+\cdots +m_{\mu}=1}
\prod_{\xi=1}^{\mu}\nu_{f_\xi,i_{\xi},m_{\xi}}
=\frac{1}{2} \, \prod_{\xi:i_{\xi}=0} (-f_{\xi})
\, \prod_{\xi:i_{\xi}=1} 2.$$
If $n=2,\, l=2$, then
$$\sum_{m_{1}+\cdots +m_{\mu}=0}
 \prod_{\xi=1}^{\mu}\nu_{f_\xi,i_{\xi},m_{\xi}}
=\prod_{\xi:i_{\xi}=0} (-f_{\xi})
\, \prod_{\xi:i_{\xi}=1} 2
\, \prod_{\xi:i_{\xi}=2} (-2f_\xi).$$
This verifies the formula.
\end{proof-5-1}

\section{Integrality of $\lambda_{n}$}

Rozansky predicts that there is an integer constant $C$ such that
$C\lambda_{n}\in\Bbb Z$.  The  predicted constant $C$ is
very large.  Conjecture 1.1 says that $C=n!$ is sufficient.  In this
section, we prove the following theorem
which is the case $n=2$ of Conjecture 1.1.

\begin{theorem}
Let $M$ be a homology 3-sphere.  Then $\lambda_{2}(M) \in 3\Bbb Z.$
\end{theorem}

If we write $\tau_{r}(M)$ as follows:
$$
\tau_r(M)=a_{r,0}+a_{r,1}(q-1)+\cdots+a_{r,n}(q-1)^n+\cdots
            +a_{r,N}(q-1)^{N}
$$
for some $a_{r,n}\in {\Bbb Z}$.
Although $a_{r,n}$ is not well-defined,
$\overline{a_{r,n}}\in {\Bbb Z}/r{\Bbb Z}$ is well-defined for
$0\leq n \leq r-2$.

\begin{corollary}
For $n=1,2$,
$\overline{a_{r,n}}\in {\Bbb Z}/r{\Bbb Z}$ stabilizes
for large $r$.
\end{corollary}

See the remark after Lemma 3.1.

It seems to be an interesting question as of when this magic stabilization
starts. We suspect that this may happen from the very first term.

Theorem 7.1 will follow from Theorem 5.1 and the following theorem, where
by doubling a component of a link $L$, we mean to split that component of $L$
into two parallel copies with zero linking number.

\begin{theorem}
Suppose $L'$ is obtained from an ASL link $L$ by doubling $m$
components $(m\geq 1)$ of $L$.  Then $\phi_{2}(L') \in 2^{m}\Bbb Z$.
\end{theorem}

To prove this theorem, we will use the {\it universal invariant} (or the
{\it colored Jones polynomial}) for ASL's. We will follow
the exposition of [O1] (see also [RT, KM, MM]).

As before, let $r$ be an odd prime and $q=e^{\frac{2 \pi\sqrt{-1}}{r}}$.  We assume
$$q^{\frac12}=-q^{\frac{r+1}{2}}.$$
Let
$U_{q}$ be the associative algebra over $\Bbb C$ with generators
$E,F, K^{\pm}$ subject to the following relations:
$$ \begin{array}{c}
K K^{-1}=K^{-1} K=1\\
KE=q EK, \;\; KF={q}^{-1}FK\\
EF-FE=\frac{K-K^{-1}}{{q}^{\frac12}-{q}^{-\frac12}}\\
E^{r}=F^{r}=0.
\end{array}
$$
Note that we do not impose a relation $K^{2r}=1$ on $U_{q}$.
The algebra $U_{q}$ becomes a Hopf algebra with comultiplication
$\Delta: U_{q}\rightarrow U_{q}\otimes U_{q}$, antipode
$S:U_{q}\rightarrow U_{q}$ and counit $\epsilon :
U_{q}\rightarrow\Bbb C$ by :
$$\begin{array}{c}
 \Delta (E)=E\otimes 1 + K\otimes E\\
\Delta (F)=F\otimes K^{-1} + 1\otimes F\\
\Delta (K^{\pm})=K^{\pm} \otimes K^{\pm}\\
S(E)=-K^{-1}E,\;\;S(F)=-FK,\;\; S(K^{\pm})=K^{\mp}\\
\epsilon (E)=\epsilon (F)=0,\;\; \epsilon (K^{\pm})=1.
\end{array}
$$
Let $I$ be the vector subspace of $U_{q}$ generated by $ab-ba$ for
any $a,b \in U_{q}$.  Then the universal invariant of a 0-framed
ASL $L$, which will be denoted by $\phi(L)$, takes value in
$( U_{q}/I )^{\otimes (\# L)}$.  By Proposition 5.6 of [O1], the
universal invariant $\phi(L)$ of $L$ has
a form obtained by taking sum, product and tensor product in
$E,F, K^{\pm}$ with coefficients in ${\Bbb Z}[q]$.

Let $\rho : U_{q} \rightarrow\text{End}({\Bbb C}^{2})$ be
the representation given by
$$
\rho (E)=
\left( \begin{array}{cc}
0 & 1\\
0 & 0
\end{array} \right) ,\;\;\;\;
\rho (F)=
\left( \begin{array}{cc}
0 & 0\\
1 & 0
\end{array} \right) ,\;\;\;\;
\rho (K^{\pm})=
\left( \begin{array}{cc}
 {q}^{\pm 1/2} & 0\\
0 & {q}^{\mp 1/2}
\end{array} \right) . $$
And let $\chi_{V} : U_{q}\rightarrow {\Bbb C}$ be the character of
$\rho$.  Note that $\chi_{V}$ descends to a map $U_{q}/I
\rightarrow\Bbb C$.
Then by Lemma 6.1 of [O1], the Jones polynomial at $q^{-1}$ for $L$
is given by
$$V(L; {q}^{-1})=\frac{(-1)^{\# L}}{q^{\frac12}+q^{-\frac12}}\,
\chi_{V}^{\otimes (\# L)} (\phi(L)).$$

We collect some facts useful in the proof of Theorem 7.2  into the following
lemma.

\begin{lemma}
Let $L'$ be as in Theorem 7.2 and assume that the doubled $m$ components
of $L$ are the first $m$ components.
\begin{enumerate}
\item $(-q^{\frac12} -q^{-\frac12})^{\# L'} \Phi(L';q^{-1})=
(\chi_{V} + (q^{\frac12} +q^{-\frac12})\epsilon)^{\otimes (\# L')}
(\phi(L'))$;
\item $\phi(L')=(\Delta \otimes \cdots \otimes \Delta \otimes id
\otimes \cdots \otimes id )(\phi(L))$ (there are $m$ $\Delta$'s here);
\item $\phi(L)$ is a ${\Bbb Z}[q]$-linear combination of
$$(q -1)^{2 j_{1}+\sum_{\nu=1}^{\mu} j_{\nu} }  K^{i_{1}}E^{j_{1}}F^{j_{1}}
\otimes K^{i_{2}}E^{j_{2}}F^{j_{2}}\otimes
\cdots \otimes K^{i_{\mu}}E^{j_{\mu}}F^{j_{\mu}} $$
such that {\em (a)} $i_{\nu}+j_{\nu}$ is even; and {\em (b)} if $j_{1}=0$, then
$$K^{i_{1}}\otimes  K^{i_{2}}E^{j_{2}}F^{j_{2}}\otimes \cdots \otimes
K^{i_{\mu}}E^{i_{\mu}}F^{i_{\mu}}$$
and
$$K^{-i_{1}}\otimes
K^{i_{2}}E^{j_{2}}F^{j_{2}}\otimes \cdots \otimes
K^{i_{\mu}}E^{i_{\mu}}F^{i_{\mu}}$$
appear in pairs with the same coefficient.
\end{enumerate}
\end{lemma}

The facts in (1), (2) and part (a) of (3) follow from Lemmas
6.2, 6.3 and 6.5  of [O1].  Part (b) of (3) follows from the
proof of Lemma 6.3 of [O1].

\begin{proof-7-2}
Theorem 7.2 follows from a careful study of the proof of Proposition 3.4 of [O1] (see
Lemma 5.2). We will examine how the universal invariant changes when we double one
component of
$L$. Note that
$\phi_{2}(L')$ is the  first possibly non-vanishing coefficient in
the Taylor expansion of $\Phi(L';t)$ at $t=1$ by Lemma 5.2.  Since
$$(-q^{\frac12} -q^{-\frac12})^{\# L'} \Phi(L';q^{-1})=
(\chi_{V} + (q^{\frac12} +q^{-\frac12})\epsilon)^{\otimes (\# L')}
(\phi(L')),$$
we need to
show that the first possibly non-vanishing coefficient in
the expansion of $(\chi_{V} + (q^{\frac12} +q^{-\frac12})\epsilon)^{\otimes (\# L')}
(\phi(L'))$ is
always a multiple of 2 when some component of $L$ is doubled.  As
observed in [O1], we can treat $q$ as an indeterminate up to a certain order
of the expansion.
Therefore, according to the argument which proves Proposition 3.4 of [O1], it is
sufficient to show that after factoring out
$(q -1)^{2}$ when $\nu \neq 1$ or $(q-1)^{4}$
when $\nu=1$ from each term
$$\begin{cases}
(q -1)^{j_{\nu}}
(\chi_{V} +(q^{\frac12} +q^{-\frac12})\epsilon )^{\otimes 2}
\Delta (K^{i_{\nu}}E^{j_{\nu}}F^{j_{\nu}}) & \text{for $\nu\neq1$}\\
 (q -1)^{2 j_{1}}
(\chi_{V} +(q^{\frac12} +q^{-\frac12})\epsilon )^{\otimes 2}
\Delta (K^{i_{1}}E^{j_{1}}F^{j_{1}}) & \text{for $\nu=1$},
\end{cases}
$$
respectively, we always have a factor 2 in the constant coefficient.

First we consider $\nu \neq 1$.
There
are three cases: $j_{\nu}=0,1,2.$

If $j_{\nu}=0$, then
$$ (\chi_{V} +(q^{\frac12} +q^{-\frac12})\epsilon )^{\otimes 2}
\Delta (K^{i_{\nu}}) =O( (q-1)^{4};\Bbb Z).$$
Since this is a higher order term, its contribution to $\phi_{2}(L')$ is 0.

If $j_{\nu}=1$, then
$$(q -1) (\chi_{V} +(q^{\frac12} +q^{-\frac12})\epsilon )^{\otimes 2}
\Delta (K^{i_{\nu}}EF)=O((q-1)^{3};\Bbb Z).$$
Again this is a higher order term.

If $j_{\nu}=2,$ then
$$(q -1)^{2} (\chi_{V} +(q^{\frac12} +q^{-\frac12})\epsilon )^{\otimes 2}
\Delta (K^{i_{\nu}}E^{2}F^{2})=
(q-1)^{2} (2+ q +q^{-1}) q^{i_{\nu}}.$$
So we have the desired factor 2.

Now consider $\nu=1$.   There
are also three cases: $j_{\nu}=0,1,2.$

If $j_{1}=0$, then
$$
\begin{aligned}
(\chi_{V} +&(q^{\frac12} +q^{-\frac12})\epsilon )^{\otimes 2}
\Delta (K^{i_{1}}) \\
&=\left[\left(\frac{i_{1}}{2}\right)^{2}-\left(\frac{r+1}{2}\right)^{2}
\right](q-1)^{4}+O((q-1)^5;{\Bbb Z})
\end{aligned}.$$
As $K^{i_{1}}$ and $K^{-i_{1}}$ appear in pairs with the same
coefficient, so we have a factor 2.

If $j_{1}=1,$ then
$$
\begin{aligned}(q -1)^{2} &(\chi_{V} +(q^{\frac12} +q^{-\frac12})\epsilon
)^{\otimes 2}
\Delta (K^{i_{1}}EF)\\
&=2(q-1)^{2} q^{\frac{i_{1}}2} (q^{\frac{i_{1}-1}{2}} + q^{-\frac{i_{1}-1}{2}}
- q^{\frac{r+1}{2}}-q^{-\frac{r+1}{2}}).
\end{aligned}$$

If $j_{1}=2,$ then
$$
\begin{aligned}
(q-1)^{4} &(\chi_{V} +(q^{\frac12} +q^{-\frac12})\epsilon
)^{\otimes 2}
\Delta (K^{i_{1}}E^{2}F^{2})\\
&=
(q-1)^{4} (2+ q +q^{-1}) q^{i_{1}}.
\end{aligned}$$
In both cases, we have the desired factor 2.
This completes the proof of Theorem 7.2.
\end{proof-7-2}

\begin{proof-7-1} Combine (5.2), Corollary 5.1 (2) and Theorem 7.2,
we get Theorem 7.1.
\end{proof-7-1}

\section{Some problems}

Here are several problems raised during the course of this work.

\noindent{\bf Problem 1:} The motivation of [O2] is to give a definition of "finite
type invariants of homology 3-spheres" which will include
$\lambda_{n}$.  To the best of our knowledge, it is still unknown whether
$\lambda_n$'s are of finite type. See [R1,2] for the physical evidence
indicating that
$\lambda_n$'s are likely to be of \lq\lq finite type".
We will address this
problem in a sequel to this paper.

\noindent{\bf Problem 2:} As $\lambda_{1}$ is essentially the Casson
invariant and
$\lambda_{2}$ is an integer, it is natural to ask if there is an analogous
interpretation for $\lambda_{2}$ as the Casson invariant; or whether
$\lambda_2/3$  counts algebraically any geometrical or topological objects
related with the manifolds in question?

\noindent{\bf Problem 3:} The series $\sum \lambda_{n}(t-1)^{n}$ is
only a formal series by definition.  It seems interesting to know
if this series actually converges to an analytic function for a
fixed homology 3-sphere.  See [L] for some examples where we do have
the convergence. This is related to Problem 1 as follows.
As conjectured in [LW] and proved in [B], for any Vassiliev invariant
$v_{k}$ of order $k$, there is a constant $C$ such that for any knot
of $n$ crossings, $|v_{k}(K)|\leq Cn^{k}$.  If $\lambda_{n}$'s are
indeed finite type invariants, it seems possible that there are
analogous inequalities for $n!\cdot \lambda_{n}$. With some control
over the constants, this will imply
that $\sum \lambda_{n} (t-1)^{n}$ converges to an analytic
function for a fixed homology 3-sphere.

\noindent{\bf Problem 4:} The
Reshetikhin-Turaev invariant $\tau_{r}(M)$ for an arbitrary 3-manifold $M$
normalized as in [KM] has the property $\tau_{r}(M)
\in{\Bbb Z}_r[q]$.  Does their Fermat limit exist?


\begin{thebibliography}{99}

\bibitem[AM]{AM}S. Akbulut and J. McCarthy, Casson's invariant for
oriented homology $3$-spheres: An exposition, Math. Notes, vol. 36,  Princeton
University Press, 1990.
\bibitem[B]{B}D. Bar-Natan, {\em Polynomial invariants are polynomial,} Math.
Res. Letter, {\bf 2}(1995).
\bibitem[Bo]{Bo}H. Boden, {\em Note on a conjecture of Lin and Wang,}
preprint, April 1996.
\bibitem[J]{J}V. Jones, {\em Hecke algebra representations
of braid groups and link polynomials,}
Ann. of Math., {\bf 186}(1987).
\bibitem[Ha]{Ha}N. Habegger, {\em Finite type 3-manifold invariants:
a proof of a conjecture of Garoufalidis,} preprint, July 1995.
\bibitem[HW]{HW}G. Hardy and E. Wright, An introduction to
the theory of numbers, Oxford University press, Fifth Edition, 1979.
\bibitem[Ho]{Ho}J. Hoste, {\em A formula for Casson's invariant,}
Tran. AMS, {\bf 297}(1986).
\bibitem[KM]{KM}R. Kirby and P. Melvin, {\em The 3-manifold
invariants of Witten and Reshetikhin-Turaev for sl(2,C),} Invent.
Math., {\bf 105}(1991).
\bibitem[G]{G}S. Garoufalidis, {\em On finite type 3-manifold
invariants I,} preprint, February 1995.
\bibitem[GLe1]{GLe1}S. Garoufalidis and J. Levine,
{\em On finite type 3-manifold invariants II,} preprint, June 1995.
\bibitem[GLe2]{GLe2}S. Garoufalidis and J. Levine,
{\em On finite type 3-manifold invariants IV: comparison
of definitions,} preprint, September 1995.
\bibitem[GLi]{GLi}M. Greenwood and X.-S. Lin,
{\em On Vassiliev knot invariants induced from finite
type 3-manifold invariants,} preprint, May 1995.
\bibitem[GO]{GO}S. Garoufalidis and T. Ohtsuki,  {\em On finite type 3-manifold
invariants III: manifold weight systems, } preprint, August 1995.
\bibitem[L]{L}R. Lawrence, {\em Asymptotic expansion of
Witten-Reshetikhin-Turaev invariants for some simple 3-manifolds}, to appear
in Jour. Math. Phys..
\bibitem[LW]{LW}X.-S. Lin and Z. Wang, {\em Integral geometry of plane
curves and knot invariants}, to appear in Jour. of Diff. Geom..
\bibitem[M]{M}H. Morton, {\em  The coloured Jones function and Alexander
polynomial for torus knots,} Math. Proc. Camb. Phil. Soc., {\bf 117}(1995).
\bibitem[MM]{MM}P. Melvin and H. Morton, {\em  The coloured Jones function,}
Comm. Math. Phys., {\bf 169}(1995).
\bibitem[Mu]{Mu}H. Murakami, {\em Quantum SU(2)-invariants dominate
Casson's SU(2)-invariant,} Math. Proc. Camb. Phil. Soc., {\bf 115}(1994).
\bibitem[O1]{O1}T. Ohtsuki, {\em A polynomial invariant of integral
homology 3-spheres,} Math. Proc. Camb. Phil. Soc., {\bf 117}(1995).
\bibitem[O2]{O2}T. Ohtsuki, {\em Finite type invariants of integral
homology 3-spheres,} preprint, 1994.
\bibitem[R1]{R1}L. Rozansky, {\em The trivial connection
contributions to Witten's invariants of rational homology
spheres,}
preprint, 1995.
\bibitem[R2]{R2}L. Rozansky, {\em Witten's invariants of rational homology
spheres at prime values of K and trivial connection contribution,}
preprint, 1995.
\bibitem[RT]{RT}N.Y. Reshetikhin and V.G. Turaev, {\em Invariants of
 3-manifolds via link polynomials and quantum groups,} Invent.
Math., {\bf 103}(1991).
\end{thebibliography}
\end{document}